\documentclass[final]{dmtcs-episciences}

\usepackage[utf8]{inputenc}
\usepackage[T1]{fontenc}
\usepackage{amsmath,amssymb}
\let\dmtcsqed\qed
\let\qed\relax
\usepackage{amsthm}
\let\qed\dmtcsqed
\usepackage{mathtools}
\usepackage[round]{natbib}
\usepackage{etoolbox}
\usepackage{tikz}
\usetikzlibrary{positioning}

\newtheorem{theorem}{Theorem}[section]

\newtheorem{lemma}[theorem]{Lemma}

\newtheorem{fact}{Fact}
\newtheorem{observation}[theorem]{Observation}
\theoremstyle{definition}

\newtheorem{claim}{Claim}

\newtheorem{case}{Case}[section]

\newtheorem{construction}[theorem]{Construction}
\AtBeginEnvironment{proof}{\setcounter{case}{0}}

\title[Nonhamiltonian locally linear graphs]{Extremal problems about the order and size of nonhamiltonian locally linear graphs}
\author[Feng Liu and Leilei Zhang]{Feng Liu\affiliationmark{1}\thanks{Feng Liu is supported by National Key R\&D Program of China under Grant No. 2022YFA1006400 and National Natural Science Foundation of China under Grant No. 12571376.}
  \and Leilei Zhang\affiliationmark{2}\thanks{Leilei Zhang's work was supported by the National Natural Science Foundation of China (Grant Nos. 12601664, 12571365), JSPS KAKENHI (Grant No.~25KF0036),  The Natural Science Foundation of Hubei Province (Grant No.~2025AFB309), the China Postdoctoral Science Foundation Grant (No.~2025M773113), and the Fundamental Research Funds for the Central Universities, Central China Normal University (Grant No. CCNU24XJ026).}}
\affiliation{%
  School of Mathematical Sciences, Shanghai Jiao Tong University, Shanghai, China\\
  School of Mathematics and Statistics, and Hubei Key Lab-Math. Sci., Central China Normal University, Wuhan, China}

\begin{document}
\publicationdata{vol. 28:3}{2026}{13}{10.46298/dmtcs.18029}{2026-04-16; 2026-04-16; 2026-09-10}{2026-09-12}
\maketitle
\begin{abstract}
		The relation between local structure and global cycle properties is a classical topic in graph theory. A graph $G$ is locally linear if $G[N(v)]$ is a path for every $v\in V(G)$. It is locally Hamiltonian or locally traceable if every vertex neighborhood induces a Hamiltonian or traceable graph, respectively. Earlier work by Pareek and Skupie\'{n}, Skupie\'{n}, Davies and Thomassen, Asratian and Oksimets, and de Wet and van Aardt studied extremal questions for these related graph classes. We prove that the minimum order of a nonhamiltonian locally linear graph is $12$ and that, for every integer $n\geq 12$, the minimum size of such a graph of order $n$ is $2n$. We also prove that every nontraceable locally linear graph of order $n$ has at least $2n+3$ edges.
	\end{abstract}
\noindent{\small\textbf{AMS Subject Classification:} 05C35, 05C45, 05C69.\par}

	\section{Introduction}
	All graphs considered in this paper are finite, undirected, simple, connected, and of order at least $3$. We use the notation and terminology of \citet{B-M}. For a graph $G$, its vertex and edge sets are denoted by $V(G)$ and $E(G)$, respectively. Its \emph{order} is $|V(G)|$, and its \emph{size} is $|E(G)|$. For $S\subseteq V(G)$, we denote by $G[S]$ the subgraph induced by $S$. A \emph{Hamilton cycle} is a cycle containing every vertex, and a graph containing such a cycle is \emph{Hamiltonian}. A \emph{Hamilton path} is a path containing every vertex, and a graph containing such a path is \emph{traceable}.

	Let $v$ be a vertex in a graph $G$. The set of all neighbors of $v$, that is, the vertices adjacent to $v$, is denoted by $N_G(v)$; when no ambiguity arises, we write $N(v)$ for simplicity. A graph \(G\) is called \emph{locally} \(\mathcal{P}\) with respect to a graph property \(\mathcal{P}\) if for each vertex \(v \in V(G)\), the subgraph induced by \(N(v)\) possesses property \(\mathcal{P}\). For example, $G$ is locally Hamiltonian if each neighborhood subgraph $G[N(v)]$ contains a Hamilton cycle. This notion was originally introduced by \citet{S1}. For further advancements regarding locally Hamiltonian and locally traceable graphs, we direct the interested reader to~\citet{JDMS,JWS,5,8}. A graph $G$ is \emph{locally linear} if $G[N(v)]$ is a path for every $v\in V(G)$.

	As pointed out by \citet*{ltz}, deducing global properties from local neighborhood structures is a classical theme in graph theory. A graph is \emph{locally forest} if every vertex neighborhood induces a forest. \Citet{moon} initiated the study of extremal problems for this class by determining its maximum size. More recent work uses locally forest graphs to study sufficient size conditions for the existence of forest vertex cuts~\citep{CRR2025} and independent minimum vertex cuts~\citep{ltz}. Since a path is a forest and is also traceable, locally linear graphs provide a natural setting in which to study how these two local conditions influence global paths and cycles.

	By definition, locally linear graphs are locally traceable, and no locally linear graph is locally Hamiltonian. Graphs with this local condition have been of interest for decades. \Citet{P-P} investigated the topological properties of such graphs. Subsequently, \citet{L-M} characterized all locally linear graphs that are free of 3-suns, and also considered cases where $G[N(v)] = P_k$ for $k \leq 4$. For broader results on neighborhood-induced substructures, see~\citet{2,3,4,6,7,CRR2025,WY2025}.

	A fundamental problem in this area is to determine the minimum order of a locally $\mathcal{P}$ graph that does not itself have property $\mathcal{P}$. \Citet{P-S2} determined the minimum order of a nonhamiltonian locally Hamiltonian graph. \Citet{JWS} proved that the minimum order of a nontraceable locally traceable graph is $10$, whereas the corresponding minimum for locally Hamiltonian graphs is $14$. These results answer two questions posed in~\citet{P-S2}. Moreover, Theorem~2.4 of~\citet{JWS} implies that the minimum order of a nonhamiltonian locally traceable graph is $7$. We determine the corresponding minimum for locally linear graphs.

	\begin{theorem}\label{MTHM1}
		The smallest nonhamiltonian locally linear graph has order $12$.
	\end{theorem}

	The lower bound of $3n-6$ edges for locally Hamiltonian graphs of order $n$ goes back to \citet{S3}. \Citet{9} completed the details of the earlier proof and extended the result to graphs with multiple edges. Their theorem also characterizes equality by the existence of an embedding that triangulates the sphere. \Citet{A-O} proved that the minimum size of a locally traceable graph of order $n$ is $2n-3$. \Citet{JWS} asked for the minimum size of a nontraceable locally traceable graph of order $n$. We study the corresponding size problems within the class of locally linear graphs and prove the following results.

	\begin{theorem}\label{MTHM2}
		For every integer $n\geq 12$, the minimum size of a nonhamiltonian locally linear graph of order $n$ is $2n$.
	\end{theorem}

	\begin{theorem}\label{MTHM3}
		The number of edges in any nontraceable locally linear graph with $n$ vertices is at least $2n + 3$.
	\end{theorem}

	\section{Preliminaries}\label{section-pre}

	In this section, we introduce the necessary notation and concepts that will be used throughout the paper.
	For clarity and ease of reference, the basic symbols are summarized in Table~\ref{tab:notation}.

	\begin{table}[htbp]
\centering
\caption{Notation used in the paper.}\label{tab:notation}
		\begin{tabular}{ll}
			\hline
			Symbol & Meaning \\
			\hline
			$N_G[v]$ & closed neighborhood of $v$: $N_G(v) \cup \{v\}$ \\
			$d_G(v)$ & degree of vertex $v$ in graph $G$, \textit{i.e.}, $|N_G(v)|$ \\
			$N_G(S)$ & $\{x\in V(G)\setminus S:x\text{ has a neighbor in }S\}$ \\
			$d_G(u,v)$ & distance between vertices $u$ and $v$ in $G$ \\
			$\Delta(G)$ & maximum degree of $G$ \\
			$\delta(G)$ & minimum degree of $G$ \\
			$C_n$ & cycle graph on $n$ vertices \\
			$K_n$ & complete graph on $n$ vertices \\
			$[k]$ & set of integers $\{1,2,\ldots,k\}$ \\
			$[a,b]$ & set of integers $\{c \in \mathbb{Z} \mid a \le c \le b\}$ \\
			\hline
		\end{tabular}
	\end{table}

	For a path $P$ in $G$ and a vertex $x\in V(G)$, write
	\begin{equation}
	N_P(x)=N_G(x)\cap V(P),\qquad N_P(S)=N_G(S)\cap V(P)
	\end{equation}
	for $S\subseteq V(G)$. In particular, when $P=v_1\cdots v_t$ is induced,
	\begin{equation}
	N_P(v_i)=
	\begin{cases}
		\{v_2\},&i=1,\\
		\{v_{i-1},v_{i+1}\},&2\leq i\leq t-1,\\
		\{v_{t-1}\},&i=t.
	\end{cases}
	\end{equation}
	These formulas are used only for $t\geq2$. For a graph $H$, a vertex $x\in V(H)$ and a nonempty set $S\subseteq V(H)$, define
	\begin{equation}
	d_H(x,S)=\min\{d_H(x,y):y\in S\}.
	\end{equation}
	Distances between vertices in different components are taken to be infinite.

	For clarity, graph equality in our discussion refers to isomorphism; that is, $G_1 = G_2$ means that $G_1$ and $G_2$ are isomorphic graphs.

	We next introduce several additional definitions and concepts that are relevant to our discussion.
	A graph $G$ is said to be \emph{fully cycle extendable} if each vertex lies in a 3-cycle and every nonhamiltonian cycle in $G$ can be extended to a larger cycle by the addition of exactly one new vertex.
	A \emph{diamond} is the graph obtained by removing one edge from the complete graph on four vertices, $K_4$.
	The independence number of a graph $G$ is denoted by $\alpha(G)$.
	For vertices $u, v \in V(G)$, we write $u \sim v$ to indicate that $u$ and $v$ are adjacent, and $u \nsim v$ when they are not adjacent.
	\section{Smallest Nonhamiltonian Locally Linear Graphs}

\begin{figure}[htbp]
		\centering
		\begin{minipage}{0.31\textwidth}
			\centering
			\begin{tikzpicture}[scale=0.95]
				\coordinate (0) at (0, 0);
				\coordinate (1) at (0, 2);
				\coordinate (2) at ({sqrt(3)}, 1);
				\coordinate (3) at ({sqrt(3)}, -1);
				\coordinate (4) at (0, -2);
				\coordinate (5) at (-{sqrt(3)}, -1);
				\coordinate (6) at (-{sqrt(3)}, 1);

				\foreach \i in {0,1,2,3,4,5,6} {
					\fill[black] (\i) circle (3pt);
				}

				\foreach \i in {1,3,5} {
					\draw[line width=0.5mm] (0) -- (\i);
				}

				\draw[line width=0.5mm] (1) -- (2);
				\draw[line width=0.5mm] (2) -- (3);
				\draw[line width=0.5mm] (3) -- (4);
				\draw[line width=0.5mm] (4) -- (5);
				\draw[line width=0.5mm] (5) -- (6);
				\draw[line width=0.5mm] (6) -- (1);
				\draw[line width=0.5mm] (1) -- (3);
				\draw[line width=0.5mm] (1) -- (5);
				\draw[line width=0.5mm] (5) -- (3);
			\end{tikzpicture}
		\end{minipage} \hfill
		\begin{minipage}{0.31\textwidth}
			\centering
			\begin{tikzpicture}[scale=0.95]
				\coordinate (0) at (0, 0);
				\coordinate (1) at (0, 2);
				\coordinate (2) at ({sqrt(2)},{sqrt(2)});
				\coordinate (3) at (2, 0);
				\coordinate (4) at ({sqrt(2)},{-sqrt(2)});
				\coordinate (5) at (0, -2);
				\coordinate (6) at ({-sqrt(2)},{-sqrt(2)});
				\coordinate (7) at (-2, 0);
				\coordinate (8) at ({-sqrt(2)},{sqrt(2)});

				\foreach \i in {0,1,...,8} {
					\fill[black] (\i) circle (3pt);
				}

				\draw[line width=0.5mm] (8) -- (2);
				\draw[line width=0.5mm] (0) -- (2);
				\draw[line width=0.5mm] (2) -- (4);
				\draw[line width=0.5mm] (0) -- (4);
				\draw[line width=0.5mm] (4) -- (6);
				\draw[line width=0.5mm] (0) -- (6);
				\draw[line width=0.5mm] (6) -- (8);
				\draw[line width=0.5mm] (0) -- (8);

				\foreach \i/\j in {1/2, 2/3, 3/4, 4/5, 5/6, 6/7, 7/8, 8/1} {
					\draw[line width=0.5mm] (\i) -- (\j);
				}
			\end{tikzpicture}
		\end{minipage} \hfill
		\begin{minipage}{0.31\textwidth}
			\centering
			\begin{tikzpicture}[scale=0.95]
				\def\radius{2}

				\coordinate (O) at (0, 0);

				\coordinate (1) at (0, 2);
				\fill (1) circle (2pt);

				\coordinate (2) at ({2*sin(36)},{2*cos(36)});
				\coordinate (3) at ({2*sin(72)},{2*cos(72)});
				\coordinate (4) at ({2*sin(108)},{2*cos(108)});
				\coordinate (5) at ({2*sin(144)},{2*cos(144)});
				\coordinate (6) at ({2*sin(180)},{2*cos(180)});
				\coordinate (7) at ({2*sin(216)},{2*cos(216)});
				\coordinate (8) at ({2*sin(252)},{2*cos(252)});
				\coordinate (9) at ({2*sin(288)},{2*cos(288)});
				\coordinate (10) at ({2*sin(324)},{2*cos(324)});

				\fill[black] (2) circle (3pt);
				\fill[black] (3) circle (3pt);
				\fill[black] (4) circle (3pt);
				\fill[black] (5) circle (3pt);
				\fill[black] (6) circle (3pt);
				\fill[black] (7) circle (3pt);
				\fill[black] (8) circle (3pt);
				\fill[black] (9) circle (3pt);
				\fill[black] (1) circle (3pt);
				\fill[black] (10) circle (3pt);

				\fill[black] (O) circle (3pt);

				\foreach \i in {1,3,5,7,9} {
					\draw[line width=0.5mm] (O) -- (\i);
				}

				\foreach \i/\j in {1/2, 2/3, 3/4, 4/5, 5/6, 6/7, 7/8, 8/9,9/10,10/1} {
					\draw[line width=0.5mm] (\i) -- (\j);
				}

				\draw[line width=0.5mm] (1) -- (3);
				\draw[line width=0.5mm] (3) -- (5);
				\draw[line width=0.5mm] (5) -- (7);
				\draw[line width=0.5mm] (7) -- (9);
				\draw[line width=0.5mm] (9) -- (1);
			\end{tikzpicture}
		\end{minipage}
		\caption{The graphs $M_3,M_4,M_5$.}
		\label{fig:fig1}
	\end{figure}

To establish our principal theorems, we first formulate several lemmas alongside a pivotal observation.

	Let $M_3$, $M_4$, and $M_5$ be the graphs depicted in Figure~\ref{fig:fig1}. Each of these graphs is nonhamiltonian and locally traceable. None is locally linear: the neighborhood of the central vertex of $M_i$ induces $C_i$ for $i\in\{3,4,5\}$.

	\begin{lemma}[{\citealt*{SMO}}]\label{LMA1}
		Let $G$ be a locally traceable graph with maximum degree at most $5$. Then $G$ is fully cycle extendable exactly when it is not isomorphic to any of the graphs $M_3$, $M_4$, or $M_5$.
	\end{lemma}

	The following two facts can be easily verified.
	\begin{fact}\label{Fact-clique-cut}
		Let $G$ be a locally linear graph of order at least $3$. Then $G$ is $2$-connected. Furthermore, if $G$ has a vertex cut $S$ with $|S|=2$, then $S$ is a clique.
	\end{fact}
	\begin{fact}\label{Fact-two}
		Let $G$ be a locally linear graph. Then every edge of $G$ belongs to at most two triangles.
	\end{fact}

	\begin{observation}\label{OBS1}
		Let $G$ be a locally linear graph and $u \in V(G)$. Suppose that the neighborhood of $u$ induces a path $G[N(u)] = v_1 v_2 \cdots v_t$. Then the following properties hold:

		\begin{enumerate}[(1)]
			\item\label{obs:consecutive-neighbors} For any vertex $v$ such that $u \nsim v$ and $N(v) \subseteq N(u)$, there exists an index $i \in [t-1]$ with $\{v_i, v_{i+1}\} \subseteq N(v)$.

			\item\label{obs:diamond} If $x$ and $y$ are adjacent vertices neither adjacent to $u$, and if $N(\{x,y\}) \subseteq N(u)$, then there exists $i \in [t-1]$ such that the induced subgraph $G[\{v_i, v_{i+1}, x, y\}]$ forms a diamond.
		\end{enumerate}
	\end{observation}

	\begin{lemma}\label{LMA2}
		Every nonhamiltonian locally linear graph $G$ of order $n\geq3$ satisfies $\Delta(G)\leq n-5$. For every integer $n\geq12$, this bound is attained.
	\end{lemma}

\begin{proof}[of Lemma~\ref{LMA2}]
To the contrary, let $G$ be a nonhamiltonian locally linear graph of order $n$ with $\Delta(G)\geq n-4$. Clearly, $\Delta(G)\leq n-2$. Furthermore, the fact that $G$ is nonhamiltonian, together with Observation~\ref{OBS1}(\ref{obs:consecutive-neighbors}), yields $\Delta(G)\neq n-2$; similarly, Observation~\ref{OBS1}(\ref{obs:diamond}) yields $\Delta(G)\neq n-3$. Therefore, $\Delta(G)= n-4$. For convenience, assume $V(G) = \{v_i : i \in [0, n-1]\}$.
		Furthermore, let $v_0$ be a vertex of degree $n-4$, and let $P=G[N(v_0)]=v_1v_2\ldots v_{n-4}$. Put $X=\{v_{n-3},v_{n-2},v_{n-1}\}$. By Observation~\ref{OBS1} and Fact~\ref{Fact-two}, we are required to exclusively analyze the scenarios where the induced subgraph \(G[X]\) takes the form of either the path graph \(P_3\) or the complete graph \(K_3\).
		\begin{claim}\label{Lemma-Claim-consecutive}
			For any $i \in[n-3,\,n-1]$,  $v_i$ has at most two consecutive neighbors on $P$.
		\end{claim}

\begin{proof}[of Claim~\ref{Lemma-Claim-consecutive}]
Otherwise, there exists $j\in [2,n-5]$ such that $\{v_{j-1},v_j,v_{j+1}\}\subseteq N(v_i)$. Nonetheless, it may be inferred that the subgraph \(G[N(v_j)]\) encompasses a cycle \(v_{j-1}v_0v_{j+1}v_iv_{j-1}\), yielding a contradiction. This reasoning validates Claim~\ref{Lemma-Claim-consecutive}.
\end{proof}

		Suppose that $G[X]=P_3$. Set $G[X]=v_{n-3}v_{n-2}v_{n-1}$.  If $d_{G[N(v_{n-2})]}(v_{n-3},v_{n-1})\geq 3$, then there exists $i\in [n-4]$ such that $v_{n-3}v_iv_{i+1}v_{n-1}$ is a subpath of $G[N(v_{n-2})]$. This implies that $G$ is Hamiltonian, a contradiction.

		Therefore, $d_{G[N(v_{n-2})]}(v_{n-3},v_{n-1})=2$. Then there exists $i\in [n-4]$ such that $v_{n-3}v_iv_{n-1}$ is a subpath of $G[N(v_{n-2})]$. Now $N(v_i)=\{v_0\}\cup N_P(v_i)\cup\{v_{n-3},v_{n-2},v_{n-1}\}$. Since $G[N(v_i)]$ is connected, $N(\{v_{n-3},v_{n-1}\})\cap N_P(v_i)\neq \emptyset$. Choose $w\in N_P(v_i)$ adjacent to an endpoint $z$ of $v_{n-3}v_{n-2}v_{n-1}$. Replace the edge $v_iw$ of $v_0v_1\cdots v_{n-4}v_0$ by the path from $v_i$ through all three vertices of $X$, ending at $z$, and then to $w$. This is a Hamilton cycle,
		a contradiction. The same insertion works when $v_i$ is an endpoint of $P$.

		Suppose that $G[X]=K_3$. By symmetry, we may assume that $ N_P(v_{n-3})\neq \emptyset$. Since $G[N(v_{n-3})]$ is connected, there exists $i\in [n-4]$ such that $v_i\in N(v_{n-3})$ and $N(v_i)\cap \{v_{n-2},v_{n-1}\}\neq \emptyset$. By symmetry, assume that $v_i\sim v_{n-2}$. Clearly, $v_{i}\nsim v_{n-1}$. Now $N(v_i)=\{v_0\}\cup N_P(v_i)\cup\{v_{n-3},v_{n-2}\}$. Since $G[N(v_i)]$ is connected, $N(\{v_{n-3},v_{n-2}\})\cap N_P(v_i)\neq \emptyset$. Choose $w\in N_P(v_i)$ adjacent to a vertex
		$z\in\{v_{n-3},v_{n-2}\}$, and let $z'$ be the other vertex of this pair.
		Replacing the edge $v_iw$ of $v_0v_1\cdots v_{n-4}v_0$ by
		$v_iz'v_{n-1}zw$ gives a Hamilton cycle, a contradiction.
		Again, $w$ is an actual neighbor of $v_i$ on $P$, so endpoints cause no exception.

		For sharpness, let $n\geq12$ and put $k=n-5$. Define $H_n$ on
		\begin{equation}
		\{u,x_1,\ldots,x_k,a,b,c,d\}
		\end{equation}
		by taking all edges $ux_i$ ($1\leq i\leq k$), all edges
		$x_ix_{i+1}$ ($1\leq i<k$), and the additional edges
		\begin{equation}
		ab,bc,bd,ax_2,cx_6,dx_3,dx_5,bx_2,bx_3,bx_5,bx_6.
		\end{equation}
		Figure~\ref{fig:fig5} illustrates this family. We have
		$d(u)=k$, $d(b)=7$, and all other degrees are at most $5$.
		Thus $\Delta(H_n)=k=n-5$.
		The neighborhood of $u$ is the path $x_1\cdots x_k$, and the
		neighborhoods of $b$ and $d$ are the paths
		$ax_2x_3dx_5x_6c$ and $x_3bx_5$, respectively.
		The neighborhoods of $a,c,x_1,x_k$ are single edges.
		Those of $x_2,x_3,x_4,x_5,x_6$ are, respectively,
		\begin{equation}
		x_1ux_3ba,\quad x_4ux_2bd,\quad x_3ux_5,\quad
		x_4ux_6bd,\quad x_7ux_5bc.
		\end{equation}
		For $7\leq i<k$, the neighborhood of $x_i$ is $x_{i-1}ux_{i+1}$.
		Hence $H_n$ is locally linear.

		If $H_n$ had a Hamilton cycle, the degree-two vertices $a$ and $c$
		would force the edges $ab$ and $bc$ into it. The cycle would therefore
		use neither $bd$ nor any edge $bx_i$, forcing both $dx_3$ and $dx_5$.
		Similarly, the degree-two vertices $x_1$ and $x_k$ force both $ux_1$
		and $ux_k$, so the cycle must use $x_3x_4$ and $x_4x_5$.
		It would then contain $x_3dx_5x_4x_3$ as a separate cycle, a contradiction.
		Thus $H_n$ is nonhamiltonian for every $n\geq12$.
\end{proof}

\begin{figure}[htbp]
		\centering
		\begin{tikzpicture}[scale=0.85,every node/.style={font=\small},vertex/.style={circle,fill,inner sep=2pt}]
			\foreach \i in {1,...,6}{\node[vertex,label=below:$x_{\i}$] (x\i) at (\i-1,0) {};}
			\node[vertex,label=below:$x_k$] (xk) at (7,0) {};
			\node[vertex,label=below:$u$] (u) at (3,-2) {};
			\node[vertex,label=above:$a$] (a) at (1,2.7) {};
			\node[vertex,label=above:$b$] (b) at (3,2.7) {};
			\node[vertex,label=above:$c$] (c) at (5,2.7) {};
			\node[vertex,label=right:$d$] (d) at (3,1.3) {};
			\foreach \i in {1,...,5}{\pgfmathtruncatemacro{\j}{\i+1}\draw (x\i)--(x\j);}
			\draw[densely dotted] (x6)--(xk);
			\foreach \i in {1,...,6}{\draw (u)--(x\i);}
			\draw (u)--(xk);
			\draw[densely dotted] (u)--(6,0);
			\draw (a)--(b)--(c) (b)--(d) (a)--(x2) (c)--(x6) (d)--(x3) (d)--(x5);
			\foreach \i in {2,3,5,6}{\draw (b)--(x\i);}
		\end{tikzpicture}
		\caption{The nonhamiltonian locally linear family $H_n$, where $k=n-5\geq7$. The dotted segment represents the path from $x_6$ to $x_k$.}
		\label{fig:fig5}
	\end{figure}

We are now ready to present the proof of Theorem \ref{MTHM1}.

\begin{proof}[of Theorem~\ref{MTHM1}]
Let $G$ be a nonhamiltonian locally linear graph of order at most $11$. Since none of $M_3,M_4,M_5$ is locally linear, Lemma~\ref{LMA1} implies that $\Delta(G)\geq6$; indeed, a fully cycle extendable graph is Hamiltonian. Lemma~\ref{LMA2} now gives $6\leq\Delta(G)\leq n-5$. Therefore, we only need to consider the case that $G$ has order $11$ and $\Delta(G)= 6$. For convenience, we may assume that $V(G)=\{v_i:i\in [0,10]\}$. Furthermore, let $v_0$ be a vertex of degree $6$, and let $P=G[N(v_0)]=v_1v_2\ldots v_6$. Put $X=\{v_7,v_8,v_9,v_{10}\}$.

		\begin{claim}\label{Claim-consecutive}
			For any $i \in[7,\,10]$,  $v_i$ has at most two consecutive neighbors on $P$.
		\end{claim}

\begin{proof}[of Claim~\ref{Claim-consecutive}]
Otherwise, there exists some \( j \in [2,5] \) such that the set \( \{v_{j-1}, v_j, v_{j+1}\} \) is contained in the neighborhood \( N(v_i) \). Nevertheless, this implies that the subgraph \( G[N(v_j)] \) includes a cycle \( v_{j-1}v_0v_{j+1}v_iv_{j-1} \), thereby leading to a contradiction. This line of reasoning establishes the validity of Claim~\ref{Claim-consecutive}.
\end{proof}

		\begin{claim}\label{Claim-Hamiltonian}
			$G[X]$ is nonhamiltonian.
		\end{claim}

\begin{proof}[of Claim~\ref{Claim-Hamiltonian}]
Otherwise, suppose that the cycle \(v_7v_8v_9v_{10}v_7\) constitutes a Hamilton cycle within the induced subgraph \(G[X]\). First, it should be noted that \(G[X]\) is not a clique. Consequently, \(G[X]\) must contain two nonadjacent vertices. Without loss of generality, we may assume that $v_7$ and $v_9$ are nonadjacent. Suppose that $v_8\nsim v_{10}$. If $d_{G[N(v_8)]}(v_7,v_9) \geq 3$, then by symmetry, there exists $i \in [6]$ such that $v_i \sim v_7$ and $v_{i+1} \sim v_8$. However, it follows that $v_0v_1\ldots v_iv_7v_{10}v_9v_8v_{i+1}\ldots v_6v_0$ is a Hamilton cycle, a contradiction. Therefore, $d_{G[N(v_8)]}(v_7,v_9)=2$. Then there exists $i\in [6]$ such that $v_7v_iv_9$ is a subpath of $G[N(v_8)]$. Clearly, $v_{10}\nsim v_i$.  Now, $N(v_i)=\{v_0\}\cup N_P(v_i)\cup\{v_7,v_8,v_9\}$.  Since $G[N(v_i)]$ is a path, $N(\{v_7,v_9\})\cap N_P(v_i)\neq \emptyset$. In either case, we can deduce that $G$ is Hamiltonian, a contradiction.  Therefore, $v_8\sim v_{10}$. By Fact~\ref{Fact-clique-cut}, we have that $N_P(\{v_8,v_{10}\})\neq \emptyset$. By symmetry, we may assume that $N_{P}(v_8)\neq \emptyset$. Now, $v_7v_{10}v_9$ is a subpath of $G[N(v_8)]$. Then there exists $i\in [6]$ such that $v_iv_7v_{10}v_9$ is a subpath of $G[N(v_8)]$ by symmetry. Now, $N(v_i)=\{v_0\}\cup N_P(v_i)\cup\{v_7,v_8\}$. Therefore, $ N(\{v_7,v_8\})\cap N_P(v_i)\neq \emptyset$. But then $G$ contains a Hamilton cycle, contradicting our assumption.
		This proves Claim~\ref{Claim-Hamiltonian}.
\end{proof}

		\begin{claim}\label{Claim-traceable}
			$G[X]$ is nontraceable.
		\end{claim}

\begin{proof}[of Claim~\ref{Claim-traceable}]
In the alternative scenario, postulate that the vertices \(v_7, v_8, v_9, v_{10}\) form a Hamilton path within the subgraph \(G[X]\). According to Claim~\ref{Claim-Hamiltonian}, it holds that \(\min\{d_{G[X]}(v_7), d_{G[X]}(v_{10})\} = 1\). Without loss of generality, let \(d_{G[X]}(v_7) = 1\).

		Suppose that $v_{10}\sim v_8$. Now $\{v_7,v_9,v_{10}\}\subseteq N(v_8)$. If $d_{G[N(v_8)]}(v_7,\,\{v_9,v_{10}\})\geq 3$, then by symmetry, there exists $i\in[6]$ such that $v_7v_iv_{i+1}v_9$ is a subpath of $G[N(v_8)]$. However, it follows that  $v_0v_1\ldots v_iv_7v_8v_{10}v_9v_{i+1}\ldots v_6v_0$ is a Hamilton cycle, a contradiction.

		Therefore, $d_{G[N(v_8)]}(v_7,\,\{v_9,v_{10}\})=2$. Then, by symmetry,  there exists $i\in [6]$ such that $v_7v_iv_9$ is a subpath of $G[N(v_8)]$. Now $N(v_i)=\{v_0\}\cup N_P(v_i)\cup\{v_7,v_8,v_9\}$. Since $G[N(v_i)]$ is a path, $N(\{v_7,v_9\})\cap N_P(v_i)\neq \emptyset$. In either case, we can deduce that $G$ is Hamiltonian, a contradiction.

		Therefore, $v_8\nsim v_{10}$. Suppose that $v_7\sim v_1$. We assert that if $N_P(v_{10}) \neq \emptyset$, then $N_P(v_{10}) = \{v_1\}$. Indeed, if there exists $j\in[2,\,6]$ such that $v_j\sim v_{10}$, then $v_1v_7v_8v_9v_{10}v_j\ldots v_6 v_0v_{j-1}\ldots v_1$ is a Hamilton cycle, a contradiction. Suppose that $v_{10}\sim v_1$. Similarly, we have that $N_P(v_7)=\{v_1\}$. Since $G[N(v_7)]$ is connected, $v_1\sim v_8$.  Similarly, $v_1\sim v_9$.

		Now, $N(v_1) =\{v_0, v_2, v_7, v_8, v_9, v_{10}\}$. This implies  that $G[N(v_1)]$ is  not connected, a contradiction. Therefore, $v_7\nsim v_1$. By symmetry, we have that
		\begin{flalign}\label{eq-1-6}
			N(\{v_7,v_{10}\})\cap \{v_1,v_6\}=\emptyset.
		\end{flalign}

		Suppose $d_{G[N(v_8)]}(v_7,v_9)\geq 3$. Then by symmetry, there exists $i\in[6]$ such that $v_7v_iv_{i+1}v_9$ is a subpath of $G[N(v_8)]$. By $(\ref{eq-1-6})$, we have that $ 2\leq i\leq 5$. Clearly, $N(v_{10})\cap N_P(v_i)=\emptyset$ and $N_P(v_{10})\cap N(v_9)\neq \emptyset$.   Let $v_j\in N_P(v_{10})\cap N(v_9)$.
		\begin{itemize}
			\item Suppose  $i=2$. By $(\ref{eq-1-6})$, we have that $j\in [2,\,5]$. Clearly, $v_{10}\nsim v_3$.  Since $v_2\nsim v_9$, $j\in \{4,5\}$. If $j=4$, then $v_0v_1v_2v_7v_8v_3v_9v_{10}v_4v_5v_6v_0$ is a Hamilton cycle, a contradiction. Therefore, $v_{j}=v_5$. Note that there exist two Hamilton paths:
			\begin{flalign}
				\begin{cases}
					v_9v_{10}v_5v_6v_0v_4v_3v_8v_7v_2v_1;\\
					v_9v_{10}v_5v_4v_3v_8v_7v_2v_1v_0v_6.
				\end{cases}
			\end{flalign}
			Since $G$ is nonhamiltonian, $N(v_9)\cap \{v_1,v_6\}=\emptyset$.  By Claim~\ref{Claim-consecutive}, we have that $v_9\nsim v_4$. Now $N(v_9)=\{v_3,v_5,v_8,v_{10}\}$.  Since $G[N(v_9)]$ is connected, $v_5\sim v_8$.  By Claim~\ref{Claim-consecutive}, we have that $v_8\nsim v_4$. By~(\ref{eq-1-6}), we have that $v_{10}\nsim v_6$. Note that there exist two Hamilton paths
			\begin{flalign}
				\begin{cases}
					v_8v_7v_2v_1v_0v_4v_3v_9v_{10}v_5v_6;\\
					v_{10}v_5v_6v_0v_1v_2v_7v_8v_9v_3v_4.
				\end{cases}
			\end{flalign}
			Therefore, $N(\{v_8,v_{10}\})\cap \{v_4,v_6\}=\emptyset$. However, it follows that $G[N(v_5)]$ is not connected, a contradiction.

			\item Suppose  $i=3$. Since $v_9\nsim v_3$, $j\neq 3$. Clearly, $N(v_{10})\cap \{v_2,v_4\}=\emptyset$.  By~(\ref{eq-1-6}), we have that $j=5$. However, it follows that $v_0v_1v_2v_3v_7v_8v_4v_9v_{10}v_5v_6v_0$ is a Hamilton cycle, a contradiction.
			\item Suppose $i=4$. This case is similar to the case when $i=2$.

			\item Suppose $i=5$. Since $v_9\nsim v_5$, $j\neq 5$. By~(\ref{eq-1-6}), we have that $j\leq 4$. However, it follows that $v_6v_9v_{10}v_j\ldots v_1v_0v_{j+1}\ldots v_5v_7v_8v_6 $ is a Hamilton cycle, a contradiction.
		\end{itemize}
		Therefore, we have that $d_{G[N(v_8)]}(v_7,v_9)=2$. By symmetry, we have that $d_{G[N(v_9)]}(v_8,v_{10})=2$. Now let $i\in [6]$ such that $v_7v_iv_9$ is a subpath of $G[N(v_8)]$ and let $j\in [6]$ such that $v_8v_jv_{10}$ is a subpath of $G[N(v_9)]$.  Assuming $i < j$ does not reduce generality.
		Clearly, $j\geq i+2$.
		By symmetry, we only need to consider the case where $i=2$ and $j\in \{4,5\}$.

		\begin{itemize}
			\item Suppose $i=2$ and $j=4$.  By Claim~\ref{Claim-consecutive}, we have that  $N(v_3)\cap \{v_8,v_9\}=\emptyset$. Clearly, $N(v_3)\cap \{v_7,v_{10}\}=\emptyset$. By~(\ref{eq-1-6}), we have that $v_9\sim v_1$. Now, there exist two Hamilton paths:
			\begin{flalign}
				\begin{cases}
					v_8v_7v_2v_1v_9v_{10}v_4v_3v_0v_6v_5;\\
					v_{10}v_9v_1v_2v_7v_8v_4v_3v_0v_6v_5.
				\end{cases}
			\end{flalign}
			Since $G$ is nonhamiltonian, $N(\{v_8,v_{10}\})\cap \{v_3,v_5\}=\emptyset$. However, it follows that $G[N(v_4)]$ is not connected, a contradiction.

			\item Suppose that $i=2$ and $j=5$. If $v_9\sim v_1$, then there exist three Hamilton paths:
			\begin{flalign}
				\begin{cases}
					v_{10}v_9v_1v_2v_7v_8v_5v_6v_0v_3v_4;\\
					v_8v_7v_2v_1v_9v_{10}v_5v_6v_0v_3v_4;\\
					v_8v_7v_2v_1v_9v_{10}v_5v_4v_3v_0v_6.
				\end{cases}
			\end{flalign}
			Since $G$ is nonhamiltonian, $N(\{v_8,v_{10}\})\cap \{v_4,v_6\}=\emptyset$. However, it follows that $G[N(v_5)]$ is not connected, a contradiction. Therefore, $v_9\nsim v_1$. By symmetry, we have that $v_8\nsim v_6$.

			If $v_7\sim v_3$, then there exist three Hamilton paths:
			\begin{flalign}
				\begin{cases}
					v_8v_9v_{10}v_5v_6v_0v_1v_2v_7v_3v_4;\\
					v_8v_9v_{10}v_5v_4v_3v_7v_2v_1v_0v_6;\\
					v_{10}v_5v_6v_0v_1v_2v_9v_8v_7v_3v_4.
				\end{cases}
			\end{flalign}
			Since $G$ is nonhamiltonian, $N(\{v_8,v_{10}\})\cap \{v_4,v_6\}=\emptyset$. However, it follows that $G[N(v_5)]$ is not connected, a contradiction. Therefore, $v_7\nsim v_3$. By symmetry, we have that $v_{10}\nsim v_4$.  Since $G[N(v_2)]$ and $G[N(v_5)]$ are connected, $v_9\sim v_3$ and $v_8\sim v_4$. However, it follows that $v_0v_6v_5v_{10}v_9v_3v_4v_8v_7v_2v_1v_0$ is a Hamilton cycle, a contradiction.
		\end{itemize}
		Therefore, $G[X]$ is nontraceable. This proves Claim~\ref{Claim-traceable}.
\end{proof}

		\begin{claim}\label{Claim-connected}
			$G[X]$ is connected.
		\end{claim}

\begin{proof}[of Claim~\ref{Claim-connected}]
To derive a contradiction, assume that the subgraph \(G[X]\) lacks connectivity. In light of Fact~\ref{Fact-two} and Observation~\ref{OBS1}, our analysis can be confined to two scenarios: \(G[X] = K_1 \cup K_3\) and \(G[X] = K_1 \cup P_3\).
		\begin{itemize}
			\item In the case where \(G[X]\) is isomorphic to the disjoint union \(K_1 \cup K_3\), we can assume, without loss of generality, that the vertex \(v_7\) constitutes the isolated component within \(G[X]\).
			By symmetry, we may assume that $ N_P(v_8)\neq \emptyset$. Since $G[N(v_8)]$ is connected, there exists $i\in [6]$ such that $v_i\in N(v_8)$ and $N(v_i)\cap \{v_9,v_{10}\}\neq \emptyset$. By virtue of symmetry, let us assume that \(v_i\sim v_9\). Suppose that $v_i\sim v_7$. If $2\leq i\leq5$, then
			$N(v_i)=\{v_0,v_{i-1},v_{i+1},v_7,v_8,v_9\}$.
			In the path $G[N(v_i)]$, the three components initially given by
			$v_{i-1}v_0v_{i+1}$, the isolated vertex $v_7$, and the edge $v_8v_9$
			must be joined. Since $v_7$ has no neighbor in $\{v_8,v_9\}$,
			after reversing $P$ and interchanging $v_8$ and $v_9$ if necessary,
			we have $v_7\sim v_{i-1}$ and $v_9\sim v_{i+1}$.
			Replacing the segment $v_{i-1}v_iv_{i+1}$ of
			$v_0v_1\cdots v_6v_0$ by
			$v_{i-1}v_7v_iv_8v_{10}v_9v_{i+1}$ gives a Hamilton cycle, a contradiction.
			Hence $i\in\{1,6\}$. By symmetry, assume that $i=1$.
			Now, $N(v_1)=\{v_7,v_8,v_9,v_2,v_0\}$. Since $N(v_7)\cap \{v_8,v_9\}=\emptyset$, $v_2\in N(v_7) \cap N(\{v_8,v_9\})$. However, it follows that $d_{G[N(v_1)]}(v_2)\geq 3$, a contradiction. Therefore, $v_i\nsim v_7$. Now $N(v_i)=\{v_0\}\cup N_P(v_i)\cup\{v_8,v_9\}$. Since $G[N(v_i)]$ is connected, $N(\{v_8,v_9\})\cap N_P(v_i)\neq \emptyset$. By symmetry, we may assume that $v_9\sim v_{i+1}$. By Observation~\ref{OBS1}, there exists $j\in [6]$ such that $\{v_j,v_{j+1}\}\subseteq N(v_7)$. By Fact \ref{Fact-two},  we have that $j\neq i$. This means \(G\) is Hamiltonian, standing in direct contradiction to our initial assumption.

			\item In the case where \(G[X]\) is isomorphic to the disjoint union \(K_1 \cup P_3\), we can assume, without loss of generality, that the vertex \(v_7\) constitutes the isolated component within \(G[X]\), and the path \(v_8v_9v_{10}\) forms the connected component \(P_3\).
			If $d_{G[N(v_9)]}(v_8,v_{10})\geq 3$, then there exists $i\in [6]$ such that $v_8v_iv_{i+1}v_{10}$ is a subpath of $G[N(v_9)]$. By Observation~\ref{OBS1}, there exists $j\in [6]$ such that $\{v_j,v_{j+1}\}\subseteq N(v_7)$. By Fact \ref{Fact-two},  we have that $j\neq i$. This implies that $G$ is Hamiltonian, a contradiction. Therefore, $d_{G[N(v_9)]}(v_8,v_{10})= 2$. Then there exists $i\in [6]$ such that $v_8v_iv_{10}$ is a subpath of $G[N(v_9)]$. If $v_i\sim v_7$, then $i\in \{1,6\}$, since $\Delta(G)=6$. By symmetry, assume that $i=1$. Now, $N(v_1)=\{v_7,v_8,v_9,v_{10},v_2,v_0\}$. Since $N(v_7)\cap \{v_8,v_{10}\}=\emptyset$, $v_2\in N(v_7) \cap N(\{v_8,v_{10}\})$. However, it follows that $d_{G[N(v_1)]}(v_2)\geq 3$, a contradiction. Therefore, $v_i\nsim v_7$.  Now $N(v_i)=\{v_0\}\cup N_P(v_i)\cup\{v_8,v_9,v_{10}\}$. Since $G[N(v_i)]$ is connected, $N(\{v_8,v_{10}\})\cap N_P(v_i)\neq \emptyset$. By symmetry, we may assume that $v_{10}\sim v_{i+1}$. By Observation~\ref{OBS1}, there exists $j\in [6]$ such that $\{v_j,v_{j+1}\}\subseteq N(v_7)$. By Fact \ref{Fact-two},  we have that $j\neq i$. This means $G$ is Hamiltonian, which contradicts our assumption.

		\end{itemize}
		Therefore, $G[X]$ is connected. This proves Claim~\ref{Claim-connected}.
\end{proof}

		By Claims~\ref{Claim-connected} and ~\ref{Claim-traceable},  we have that $G[X]=K_{1,3}$. Without loss of generality, we may assume that $d_{G[X]}(v_7)=3$. That is, $\{v_8,v_9,v_{10}\}$ is an independent set. Since $\Delta(G)=6$, there exist $r,s\in[8,\,10]$ such that $d_{G[N(v_7)]}(v_r,v_s)=2$. By symmetry, assume that $v_8v_iv_9$ is a subpath of $G[N(v_7)]$ for some $i\in [6]$. Clearly, $v_i\nsim v_{10}$. Indeed, if $v_i\sim v_{10}$, then $d_{G[N(v_i)]}(v_7)=3$, a contradiction.
		\begin{claim}\label{Claim-claw-i}
			$i\in [2,5]$.
		\end{claim}

\begin{proof}[of Claim~\ref{Claim-claw-i}]
Otherwise, by symmetry, assume that $i=1$. Now, $N(v_1)=\{v_0,v_2,v_7,v_8,v_9\}$. Since $G[N(v_1)]$ is connected, $N(v_2)\cap \{v_8,v_9\}\neq\emptyset$. By symmetry, assume that $v_9\sim v_2$. Since $G[N(v_7)]$ is connected, there exists $j\in [6]$ such that $v_j\in N(v_7)\cap N(v_{10})$. Note that $v_2\nsim v_7 $. Therefore, $j\geq 3$.  Nonetheless, it may be deduced that the sequence \(v_0v_6\ldots v_jv_{10}v_7v_8v_1v_9v_2\ldots v_{j-1}v_0\) constitutes a Hamilton cycle, thereby yielding a contradiction. This line of reasoning establishes Claim~\ref{Claim-claw-i}.
\end{proof}

		Since $G[N(v_7)]$ is connected, choose $v_j\in N_P(v_7)\cap N(v_{10})$. Note that $i\neq j$ and $v_i\nsim v_j$. Therefore, $j\leq i-2$ or $j\geq i+2$.  Note that $N(v_i)=\{v_0\}\cup N_P(v_i)\cup\{v_7,v_8,v_9\}$. Since $G[N(v_i)]$ is connected, $N(\{v_8,v_9\})\cap N_P(v_i)\neq \emptyset$. This implies that $v_{10}$ has no  two consecutive neighbors on $P$.

		\begin{case}
			$d_{G[N(v_7)]}(v_i, v_j) = 2$.
		\end{case}

		Without loss of generality, assume $v_j \sim v_9$. By Claim~\ref{Claim-claw-i}, we have $j \in [2,5]$. Similarly, $v_8$ has no two consecutive neighbors on $P$. By symmetry, we may take $i < j$, so $j \geq i + 2$. By symmetry considerations, the analysis can be restricted to the case where \(i = 2\) and \(j\) is either 4 or 5.
		\begin{itemize}
			\item  Suppose that $i=2$ and $j=4$. Clearly, $v_9\nsim v_3$.  Since $v_8$ and $v_{10}$ have no two consecutive neighbors on $P$, we have $v_9\sim v_1$ and $v_9\sim v_5$. However, it follows that $v_0v_6v_5v_9v_1v_2v_8v_7v_{10}v_4v_3v_0$ is a Hamilton cycle, a contradiction.
			\item Suppose that $i=2$ and $j=5$.  If $v_9\sim v_3$, then Claim~\ref{Claim-consecutive} gives $v_9\nsim v_4$. The connectedness of $G[N(v_5)]$ therefore implies $v_9\sim v_6$. This implies that $v_0v_1v_2v_8v_7v_{10}v_5v_4v_3v_9v_6v_0$ is a Hamilton cycle, a contradiction.  Then $v_9\nsim v_3$.  By symmetry, we have that $v_9\nsim v_4$. This implies that $v_9\sim v_1$ and $v_9\sim v_6$. However, it follows that $v_0v_3v_4v_5v_{10}v_7v_8v_2v_1v_9v_6v_0$ is a Hamilton cycle, a contradiction.
		\end{itemize}
		\begin{case}
			$d_{G[N(v_7)]}(v_i,v_j)=3$.
		\end{case}
		Without loss of generality, let the distance between \(v_9\) and \(v_j\) in \(G[N(v_7)]\) be 2, so \(d_{G[N(v_7)]}(v_9, v_j) = 2\). Consequently, the neighborhood $G[N(v_7)]$ must be one of the two paths:
		\begin{equation}
		v_8 v_i v_9 v_{j-1} v_j v_{10} \quad \text{or} \quad v_8 v_i v_9 v_{j+1} v_j v_{10}.
		\end{equation}
		Observe that $v_i$ and $v_j$ are not adjacent ($v_i \nsim v_j$) and that $i \neq j$. This fact ensures that $v_8$ does not have two neighbors that are consecutive on the path $P$. Hence, the neighborhood of $v_9$ intersects with at least one of the vertices $v_{i-1}$ or $v_{i+1}$, \textit{i.e.}, $N(v_9) \cap \{v_{i-1}, v_{i+1}\} \neq \emptyset$. By virtue of symmetry, the analysis can be confined to the scenario where the neighborhood subgraph \(G[N(v_7)]\) is structured as the path
		\begin{equation}
		G[N(v_7)] = v_8 v_i v_9 v_{j-1} v_j v_{10}.
		\end{equation}

		\begin{itemize}
			\item Assume $i \geq j + 2$. If $v_9$ is adjacent to $v_{i-1}$, then the cycle
			\begin{equation}
			v_0 v_6 \cdots v_i v_8 v_7 v_{10} v_j \cdots v_{i-1} v_9 v_{j-1} \cdots v_1 v_0
			\end{equation}
			forms a Hamilton cycle, which contradicts our assumption. On the other hand, if $v_9$ is adjacent to $v_{i+1}$, the cycle
			\begin{equation}
			v_0 v_6 \cdots v_{i+1} v_9 v_i v_8 v_7 v_{10} v_j \cdots v_1 v_0
			\end{equation}
			is a Hamilton cycle, again leading to a contradiction.

			\item Suppose $i \leq j - 2$. If $v_9$ is adjacent to $v_{i-1}$, then the cycle
			\begin{equation}
			v_0 v_6 \cdots v_j v_{10} v_7 v_8 v_i \cdots v_{j-1} v_9 v_{i-1} \cdots v_1 v_0
			\end{equation}
			forms a Hamilton cycle, contradicting our assumptions. Hence, it must be that $v_9 \sim v_{i+1}$. Note that $i \in [2, j-2]$. From Claim~\ref{Claim-consecutive}, we deduce $i \leq j - 3$, which forces $i=2$ and $j=6$. However, this leads to the cycle
			\begin{equation}
			v_0 v_1 v_2 v_8 v_7 v_{10} v_6 v_5 v_9 v_3 v_4 v_0,
			\end{equation}
			a Hamilton cycle, again a contradiction.

		\end{itemize}
		Figure~\ref{fig:fig2} presents a connected locally linear graph with $12$ vertices. The four outer vertices of degree two would force an $8$-cycle in any Hamilton cycle, leaving the four inner vertices uncovered. Thus this graph is nonhamiltonian and the bound is sharp. This proves Theorem~\ref{MTHM1}.
\end{proof}

\begin{figure}[htbp]
		\centering
		\begin{tikzpicture}[scale=0.8,main_node/.style={circle,fill=black,minimum size=0.55em,inner sep=0pt}]
			\node[main_node] (0) at (0, 2.4) {};
			\node[main_node] (1) at (-1.6, 1.6) {};
			\node[main_node] (2) at (1.6, 1.6) {};
			\node[main_node] (3) at (-1.6, -1.6) {};
			\node[main_node] (4) at (1.6, -1.6) {};
			\node[main_node] (5) at (0, 0.8) {};
			\node[main_node] (6) at (-0.8, 0) {};
			\node[main_node] (7) at (0, -0.8) {};
			\node[main_node] (8) at (0.8, 0) {};
			\node[main_node] (9) at (-2.4, 0) {};
			\node[main_node] (10) at (0, -2.4) {};
			\node[main_node] (11) at (2.4, 0) {};

			\path[draw, thick]
			(0, 2.4) -- (-1.6, 1.6)
			(-1.6, 1.6) -- (1.6, 1.6)
			(1.6, 1.6) -- (0, 2.4)
			(-1.6, 1.6) -- (-1.6, -1.6)
			(-1.6, -1.6) -- (1.6, -1.6)
			(1.6, -1.6) -- (1.6, 1.6)
			(0, 0.8) -- (-0.8, 0)
			(-0.8, 0) -- (0, -0.8)
			(0, -0.8) -- (0.8, 0)
			(0.8, 0) -- (0, 0.8)
			(-1.6, 1.6) -- (-0.8, 0)
			(-1.6, 1.6) -- (0, 0.8)
			(0, 0.8) -- (1.6, 1.6)
			(0.8, 0) -- (1.6, 1.6)
			(0.8, 0) -- (1.6, -1.6)
			(0, -0.8) -- (1.6, -1.6)
			(-0.8, 0) -- (-1.6, -1.6)
			(-1.6, -1.6) -- (0, -0.8)
			(-2.4, 0) -- (-1.6, 1.6)
			(-2.4, 0) -- (-1.6, -1.6)
			(-1.6, -1.6) -- (0, -2.4)
			(0, -2.4) -- (1.6, -1.6)
			(1.6, 1.6) -- (2.4, 0)
			(2.4, 0) -- (1.6, -1.6);

		\end{tikzpicture}
		\caption{A connected nonhamiltonian locally linear graph of order $12$.}
		\label{fig:fig2}
	\end{figure}

\section{On the Minimum Size of Nonhamiltonian and Nontraceable Locally Linear Graphs}

	In this section, we prove Theorems~\ref{MTHM2} and~\ref{MTHM3}. \Citet{JWS} introduced an edge identification operation for locally traceable graphs. An edge $uv$ is \emph{suitable} if $G[N(u)]$ has a Hamilton path with endpoint $v$, and $G[N(v)]$ has a Hamilton path with endpoint $u$.

	\begin{construction}[{\citealt{JWS}}]
		Consider two locally traceable graphs \(G_1\) and \(G_2\), each featuring a designated edge \(u_i v_i\) for \(i = 1, 2\). A novel graph \(G\) with vertex cardinality \(|V(G_1)| + |V(G_2)| - 2\) is formed via the following procedure: Merge vertices \(u_1\) and \(u_2\) into a single vertex \(u\), and similarly coalesce vertices \(v_1\) and \(v_2\) into a single vertex \(v\). The resulting graph \(G\) retains all edges from the original graphs \(G_1\) and \(G_2\). This operation yields \(G\) as the amalgamation of \(G_1\) and \(G_2\) through the identification of their respective edges \(u_1 v_1\) and \(u_2 v_2\).
	\end{construction}

	\begin{observation}[{\citealt{JWS}}]
		In a locally traceable graph, any edge incident to a degree-two vertex $v$ is considered suitable for edge identification.
	\end{observation}

	\begin{lemma}[{\citealt*{JDMS}}]\label{MLMA3}
		Assume that $G_1$ and $G_2$ are locally traceable graphs of order at least 3, and let $G$ be constructed by identifying suitable edges of $G_1$ and $G_2$. Then the following properties are satisfied:
		\begin{enumerate}[(a)]
			\item\label{item:traceability} The graph \(G\) stays locally traceable.
			\item\label{item:planarity} In case both \(G_1\) and \(G_2\) are planar, \(G\) is also planar.
			\item\label{item:hamiltonicity} If \(G\) is Hamiltonian, then each of \(G_1\) and \(G_2\) has to be Hamiltonian too.
		\end{enumerate}

	\end{lemma}
	For an edge $e$ of a graph $G$, let $t_G(e)$ denote the number of
	triangles containing $e$, and let $t(G)$ be the number of triangles of $G$.
	In a locally linear graph of order at least three, an edge is called a
	\emph{$1$-edge} or a \emph{$2$-edge} if it belongs to exactly one or
	exactly two triangles, respectively.
	The following lemma supplies the embedding used in both proofs below.
	We do not require this embedding to be orientable.

	\begin{lemma}\label{lem:boundary-embedding}
		Let $G$ be a connected locally linear graph of order $n\geq 3$ and
		size $m$. The $1$-edges form a spanning $2$-regular subgraph of $G$.
		If this subgraph has $b$ components, then $G$ has a $2$-cell embedding
		in a connected closed surface $\Sigma$ with precisely $t(G)+b$ faces:
		one triangular face for each triangle of $G$, and one additional face
		bounded by each component of the $1$-edge subgraph.
		Writing $\Phi=\chi(\Sigma)$ for the Euler characteristic of $\Sigma$, we have $\Phi\leq 2$ and
		\begin{equation}\label{eq:boundary-count}
			3t(G)+n=2m,
			\qquad m-2n=3(b-\Phi).
		\end{equation}
	\end{lemma}

\begin{proof}
By Fact~\ref{Fact-clique-cut}, we have $\delta(G)\geq 2$.
		For a vertex $v$ and a neighbor $u$ of $v$, the number $t_G(uv)$
		equals the degree of $u$ in the path $G[N(v)]$. Thus every edge is a
		$1$-edge or a $2$-edge, and the two endpoints of $G[N(v)]$ correspond
		to exactly two $1$-edges incident with $v$. Consequently, the
		$1$-edges form a spanning $2$-regular subgraph, whose components are
		pairwise vertex-disjoint cycles $C_1,\ldots,C_b$, with $b\geq 1$.
		In particular, there are exactly $n$ $1$-edges. Counting incidences
		between edges and triangles gives
		\begin{equation}
		3t(G)=n+2(m-n)=2m-n.
		\end{equation}

		Attach a triangular disk to every triangle of $G$, identifying its
		boundary with that triangle. The resulting simplicial complex $S$
		is a compact connected surface with boundary. Indeed, an interior
		point of an edge has a disk or a half-disk neighborhood according
		as the edge is a $2$-edge or a $1$-edge. At a vertex $v$, the link
		is exactly the path $G[N(v)]$, so $v$ has a half-disk neighborhood.
		The boundary of $S$ therefore consists precisely of the cycles
		$C_1,\ldots,C_b$.

		Cap each boundary component by a disk to obtain a connected closed
		surface $\Sigma$. This gives the claimed $2$-cell embedding of $G$.
		The additional faces are the capping disks; they are distinguished
		from the original triangular faces even if their boundaries are
		triangles. Let $F(G)$ be the face set of this embedding and let $B(G)$ be the set of
		capping faces. Then $|B(G)|=b$, $|F(G)|=t(G)+b$, and
		\begin{equation}
		3|F(G)|-3|B(G)|=2m-n.
		\end{equation}
		Here $B(G)$ records the capping faces, rather than all non-triangular faces.
		Euler's formula now gives
		\begin{equation}
		\Phi=n-m+t(G)+b.
		\end{equation}
		Combining this identity with $3t(G)=2m-n$ proves
		$m-2n=3(b-\Phi)$.
		Finally, the Euler characteristic of a connected closed surface is
		at most two, whether or not the surface is orientable.
\end{proof}

\subsection{Proof of Theorem \ref{MTHM2}}

\begin{proof}
Let $G$ be a nonhamiltonian locally linear graph of order $n$ and
	size $m$, and let $b$ and $\Phi$ be as in
	Lemma~\ref{lem:boundary-embedding}. If $m<2n$, then
	\eqref{eq:boundary-count} gives
	\begin{equation}
	1\leq b<\Phi\leq 2.
	\end{equation}
	Hence $b=1$ and $\Phi=2$, so the constructed embedding is planar. Since the $1$-edges form a spanning
	$2$-regular subgraph with just one component, they form a Hamilton
	cycle of $G$, a contradiction. Therefore $m\geq 2n$.

	By Theorem \ref{MTHM1}, we have $n \geq 12$. When $n = 12$, Figure \ref{fig:fig2} shows a connected, nonhamiltonian, locally linear graph $G_{12}$ with $12$ vertices and $24$ edges. For each $i \geq 13$, define $G_i$ by attaching a $K_3$ to $G_{i-1}$ via edge identification: initially along edge $uv$, and thereafter along the edge connecting the vertices of degree two and three in the last added triangle, as illustrated in Figure \ref{fig:fig4}. Applying Lemma~\ref{MLMA3}(\ref{item:hamiltonicity}) repeatedly, each $G_i$ ($i \geq 13$) is nonhamiltonian with $i$ vertices and $2i$ edges.

	We prove inductively that every $G_i$ is locally linear. Suppose that
	$G_{k-1}$ is locally linear, where $k\geq 13$, and that $G_k$ is
	obtained by adding a vertex $x$ with neighborhood $\{y,z\}$, where
	$yz$ is the chosen edge of $G_{k-1}$. Relabel its endpoints so that
	$d_{G_{k-1}}(y)=2$, and write $N_{G_{k-1}}(y)=\{z,w\}$.
	Since $G_{k-1}$ is locally linear, $z\sim w$. The neighborhood of
	$x$ induces the path $yz$, and the neighborhood of $y$ in $G_k$
	induces the path $wzx$. In $G_{k-1}[N_{G_{k-1}}(z)]$, the vertex
	$y$ has degree one, so it is an endpoint of this induced path.
	Adding $x$ extends that path at the endpoint $y$, because $x$ is
	adjacent to $y$ and to no other vertex of $N_{G_{k-1}}(z)$.
	Thus $G_k[N_{G_k}(z)]$ is also a path. The neighborhood subgraphs
	of all other old vertices are unchanged. Therefore $G_k$ is
	locally linear, completing the induction.
	This completes the proof.
\end{proof}

\begin{figure}[htbp]
		\centering
		\begin{tikzpicture}[scale=0.8,main_node/.style={circle,draw,fill=black,inner sep=2pt}, line width=1.2pt]

			\begin{scope}
				\path[draw]
				(0, 2.4) -- (-1.6, 1.6)
				(-1.6, 1.6) -- (1.6, 1.6)
				(1.6, 1.6) -- (0, 2.4)
				(-1.6, 1.6) -- (-1.6, -1.6)
				(-1.6, -1.6) -- (1.6, -1.6)
				(1.6, -1.6) -- (1.6, 1.6)
				(0, 0.8) -- (-0.8, 0)
				(-0.8, 0) -- (0, -0.8)
				(0, -0.8) -- (0.8, 0)
				(0.8, 0) -- (0, 0.8)
				(-1.6, 1.6) -- (-0.8, 0)
				(-1.6, 1.6) -- (0, 0.8)
				(0, 0.8) -- (1.6, 1.6)
				(0.8, 0) -- (1.6, 1.6)
				(0.8, 0) -- (1.6, -1.6)
				(0, -0.8) -- (1.6, -1.6)
				(-0.8, 0) -- (-1.6, -1.6)
				(-1.6, -1.6) -- (0, -0.8)
				(-2.4, 0) -- (-1.6, 1.6)
				(-2.4, 0) -- (-1.6, -1.6)
				(-1.6, -1.6) -- (0, -2.4)
				(0, -2.4) -- (1.6, -1.6)
				(1.6, 1.6) -- (2.4, 0)
				(2.4, 0) -- (1.6, -1.6);

				\node[main_node] (0) at (0, 2.4) {};
				\node[main_node] (1) at (-1.6, 1.6) {};
				\node[main_node] (2) at (1.6, 1.6) {};
				\node[main_node] (3) at (-1.6, -1.6) {};
				\node[main_node] (4) at (1.6, -1.6) {};
				\node[main_node] (5) at (0, 0.8) {};
				\node[main_node] (6) at (-0.8, 0) {};
				\node[main_node] (7) at (0, -0.8) {};
				\node[main_node] (8) at (0.8, 0) {};
				\node[main_node] (9) at (-2.4, 0) {};
				\node[main_node] (10) at (0, -2.4) {};
				\node[main_node] (11) at (2.4, 0) {};

				\node[above=0cm of 2]{$v$};
				\node[below=0cm of 11]{$u$};
			\end{scope}

			\begin{scope}[xshift=8.5cm]
				\path[draw]
				(0, 2.4) -- (-1.6, 1.6)
				(-1.6, 1.6) -- (1.6, 1.6)
				(1.6, 1.6) -- (0, 2.4)
				(-1.6, 1.6) -- (-1.6, -1.6)
				(-1.6, -1.6) -- (1.6, -1.6)
				(1.6, -1.6) -- (1.6, 1.6)
				(0, 0.8) -- (-0.8, 0)
				(-0.8, 0) -- (0, -0.8)
				(0, -0.8) -- (0.8, 0)
				(0.8, 0) -- (0, 0.8)
				(-1.6, 1.6) -- (-0.8, 0)
				(-1.6, 1.6) -- (0, 0.8)
				(0, 0.8) -- (1.6, 1.6)
				(0.8, 0) -- (1.6, 1.6)
				(0.8, 0) -- (1.6, -1.6)
				(0, -0.8) -- (1.6, -1.6)
				(-0.8, 0) -- (-1.6, -1.6)
				(-1.6, -1.6) -- (0, -0.8)
				(-2.4, 0) -- (-1.6, 1.6)
				(-2.4, 0) -- (-1.6, -1.6)
				(-1.6, -1.6) -- (0, -2.4)
				(0, -2.4) -- (1.6, -1.6)
				(1.6, 1.6) -- (2.4, 0)
				(2.4, 1.6) -- (2.4, 0)
				(2.4, 1.6) -- (1.6, 1.6)
				(2.4, 0) -- (1.6, -1.6);

				\node[main_node] (0) at (0, 2.4) {};
				\node[main_node] (1) at (-1.6, 1.6) {};
				\node[main_node] (2) at (1.6, 1.6) {};
				\node[main_node] (3) at (-1.6, -1.6) {};
				\node[main_node] (4) at (1.6, -1.6) {};
				\node[main_node] (5) at (0, 0.8) {};
				\node[main_node] (6) at (-0.8, 0) {};
				\node[main_node] (7) at (0, -0.8) {};
				\node[main_node] (8) at (0.8, 0) {};
				\node[main_node] (9) at (-2.4, 0) {};
				\node[main_node] (10) at (0, -2.4) {};
				\node[main_node] (11) at (2.4, 0) {};
				\node[main_node] (12) at (2.4, 1.6) {};
				\node[above=0cm of 2]{$v$};
				\node[below=0cm of 11]{$u$};
			\end{scope}

			\begin{scope}[xshift=8.5cm, yshift=-8cm]
				\path[draw]
				(0, 2.4) -- (-1.6, 1.6)
				(-1.6, 1.6) -- (1.6, 1.6)
				(1.6, 1.6) -- (0, 2.4)
				(-1.6, 1.6) -- (-1.6, -1.6)
				(-1.6, -1.6) -- (1.6, -1.6)
				(1.6, -1.6) -- (1.6, 1.6)
				(0, 0.8) -- (-0.8, 0)
				(-0.8, 0) -- (0, -0.8)
				(0, -0.8) -- (0.8, 0)
				(0.8, 0) -- (0, 0.8)
				(-1.6, 1.6) -- (-0.8, 0)
				(-1.6, 1.6) -- (0, 0.8)
				(0, 0.8) -- (1.6, 1.6)
				(0.8, 0) -- (1.6, 1.6)
				(0.8, 0) -- (1.6, -1.6)
				(0, -0.8) -- (1.6, -1.6)
				(-0.8, 0) -- (-1.6, -1.6)
				(-1.6, -1.6) -- (0, -0.8)
				(-2.4, 0) -- (-1.6, 1.6)
				(-2.4, 0) -- (-1.6, -1.6)
				(-1.6, -1.6) -- (0, -2.4)
				(0, -2.4) -- (1.6, -1.6)
				(1.6, 1.6) -- (2.4, 0)
				(2.4, 1.6) -- (2.4, 0)
				(2.4, 1.6) -- (1.6, 1.6)
				(3.2, 0) -- (2.4, 1.6)
				(3.2, 0) -- (2.4, 0)
				(2.4, 0) -- (1.6, -1.6);

				\node[main_node] (0) at (0, 2.4) {};
				\node[main_node] (1) at (-1.6, 1.6) {};
				\node[main_node] (2) at (1.6, 1.6) {};
				\node[main_node] (3) at (-1.6, -1.6) {};
				\node[main_node] (4) at (1.6, -1.6) {};
				\node[main_node] (5) at (0, 0.8) {};
				\node[main_node] (6) at (-0.8, 0) {};
				\node[main_node] (7) at (0, -0.8) {};
				\node[main_node] (8) at (0.8, 0) {};
				\node[main_node] (9) at (-2.4, 0) {};
				\node[main_node] (10) at (0, -2.4) {};
				\node[main_node] (11) at (2.4, 0) {};
				\node[main_node] (12) at (2.4, 1.6) {};
				\node[main_node] (13) at (3.2, 0) {};
				\node[above=0cm of 2] {$v$};
				\node[below=0cm of 11] {$u$};
			\end{scope}

			\begin{scope}[yshift=-8cm]
				\path[draw]
				(0, 2.4) -- (-1.6, 1.6)
				(-1.6, 1.6) -- (1.6, 1.6)
				(1.6, 1.6) -- (0, 2.4)
				(-1.6, 1.6) -- (-1.6, -1.6)
				(-1.6, -1.6) -- (1.6, -1.6)
				(1.6, -1.6) -- (1.6, 1.6)
				(0, 0.8) -- (-0.8, 0)
				(-0.8, 0) -- (0, -0.8)
				(0, -0.8) -- (0.8, 0)
				(0.8, 0) -- (0, 0.8)
				(-1.6, 1.6) -- (-0.8, 0)
				(-1.6, 1.6) -- (0, 0.8)
				(0, 0.8) -- (1.6, 1.6)
				(0.8, 0) -- (1.6, 1.6)
				(0.8, 0) -- (1.6, -1.6)
				(0, -0.8) -- (1.6, -1.6)
				(-0.8, 0) -- (-1.6, -1.6)
				(-1.6, -1.6) -- (0, -0.8)
				(-2.4, 0) -- (-1.6, 1.6)
				(-2.4, 0) -- (-1.6, -1.6)
				(-1.6, -1.6) -- (0, -2.4)
				(0, -2.4) -- (1.6, -1.6)
				(1.6, 1.6) -- (2.4, 0)
				(2.4, 1.6) -- (2.4, 0)
				(2.4, 1.6) -- (1.6, 1.6)
				(3.2, 0) -- (2.4, 1.6)
				(3.2, 0) -- (2.4, 0)
				(3.2, 1.6) -- (3.2, 0)
				(3.2, 1.6) -- (2.4, 1.6)
				(2.4, 0) -- (1.6, -1.6);

				\node[main_node] (0) at (0, 2.4) {};
				\node[main_node] (1) at (-1.6, 1.6) {};
				\node[main_node] (2) at (1.6, 1.6) {};
				\node[main_node] (3) at (-1.6, -1.6) {};
				\node[main_node] (4) at (1.6, -1.6) {};
				\node[main_node] (5) at (0, 0.8) {};
				\node[main_node] (6) at (-0.8, 0) {};
				\node[main_node] (7) at (0, -0.8) {};
				\node[main_node] (8) at (0.8, 0) {};
				\node[main_node] (9) at (-2.4, 0) {};
				\node[main_node] (10) at (0, -2.4) {};
				\node[main_node] (11) at (2.4, 0) {};
				\node[main_node] (12) at (2.4, 1.6) {};
				\node[main_node] (13) at (3.2, 0) {};
				\node[main_node] (14) at (3.2, 1.6) {};
				\node[above=0cm of 2] {$v$};
				\node[below=0cm of 11] {$u$};
			\end{scope}

			\draw[->, thick, line width=1.5pt, >=stealth] (4, 0) -- (5, 0);
			\draw[->, thick, line width=1.5pt, >=stealth] (8.5, -3.5) -- (8.5, -4.5);
			\draw[->, thick, line width=1.5pt, >=stealth] (5, -8) -- (4, -8);

		\end{tikzpicture}
		\caption{Building nonhamiltonian locally linear graphs with $m = 2n$.}
		\label{fig:fig4}
	\end{figure}

\clearpage
\subsection{Proof of Theorem \ref{MTHM3}}

\begin{proof}
Let $G$ be a nontraceable locally linear graph of order $n$ and size
	$m$. Then $G$ is nonhamiltonian, so Theorem~\ref{MTHM2} gives
	$m\geq 2n$. By Lemma~\ref{lem:boundary-embedding},
	\begin{equation}
	m-2n=3(b-\Phi),
	\end{equation}
	where $b\geq 1$ is the number of components of the spanning
	$2$-regular subgraph formed by the $1$-edges, and $\Phi\leq 2$
	is an integer. Consequently, either $m=2n$ or $m\geq 2n+3$.

	Suppose that $m=2n$. Then $b=\Phi$, so $1\leq b=\Phi\leq 2$.
	If $b=1$, the $1$-edges form a Hamilton cycle, contradicting
	nontraceability. Hence $b=\Phi=2$, so the constructed embedding is planar. Let $C_1$ and $C_2$ be the two
	components of the $1$-edge subgraph. They are vertex-disjoint cycles
	whose vertex sets partition $V(G)$. Since $G$ is connected, there
	is an edge $uv$ with $u\in V(C_1)$ and $v\in V(C_2)$.
	Delete an edge incident with $u$ from $C_1$ and an edge incident
	with $v$ from $C_2$, and join the resulting paths with $uv$.
	This produces a Hamilton path of $G$, again a contradiction.
	Therefore $m\geq 2n+3$.
\end{proof}

\section{Conclusion}

	In this work, we have investigated structural properties of nonhamiltonian and nontraceable locally linear graphs.
	In particular, we have shown that the maximum possible degree of a nonhamiltonian locally linear graph of order $n$ is $n-5$ for every $n\geq12$, the minimum order of such a graph is $12$, and its minimum size is $2n$ for every $n\geq12$. Every nontraceable locally linear graph of order $n$ has at least $2n+3$ edges. We do not claim that this last bound is sharp.

	A natural direction for future research is to study the corresponding extremal problems for the maximum number of edges.
	Specifically, it would be interesting to determine the maximum size of a nonhamiltonian locally linear graph, the maximum size of a nontraceable locally linear graph, and the maximum size of a general locally linear graph for a given order.
	Investigating these dual problems could further deepen our understanding of the extremal structure of locally linear graphs.

	Although this work is primarily theoretical, the results on locally linear graphs have potential applications.
	For example, in communication networks, a locally linear neighborhood structure can model a chain of routers where each router connects linearly to a few others, helping analyze reliable message routing.
	In surveillance systems, sensors arranged in locally linear patterns can optimize coverage along corridors or pipelines.
	In distributed computing, tasks with linear dependency chains can be represented by locally linear subgraphs, facilitating scheduling and resource allocation analysis.
	These examples illustrate that the structural properties studied here may provide useful insights beyond purely theoretical interest.

	\section*{Declaration}

	\noindent\textbf{Conflict of Interest} \\
	The authors hereby declare that there are no known competing financial interests or personal relationships that could have potentially influenced the work presented in this manuscript.

	\noindent\textbf{Availability of Data} \\
	The practice of data sharing finds no applicability in the context of this investigative work, in light of the fact that no original data were either produced or subjected to analysis throughout the course of the research undertaking.

\acknowledgments

The authors would like to thank the anonymous referees for their careful reading and valuable suggestions. The authors would like to thank Professor Xingzhi Zhan for suggesting the research topic of this paper.

\nocite{1}
\bibliographystyle{abbrvnat}
\bibliography{references}
\end{document}